\documentclass{gtart_a}
\pdfoutput=1

\usepackage{pinlabel, diagrams}

\title{A note on knot Floer homology of links}
\author{Yi Ni}
\givenname{Yi}
\surname{Ni}
\address{Department of Mathematics\\Princeton University\\\newline
Princeton, NJ 08544\\USA}
\email{yni@math.princeton.edu}
\urladdr{}

\volumenumber{10}
\issuenumber{}
\publicationyear{2006}
\papernumber{19}
\lognumber{0622}
\startpage{695}
\endpage{713}

\doi{}
\MR{}
\Zbl{}

\arxivreference{math.GT/0506208}
\arxivpassword{bn2e6}  

\keyword{knot Floer homology}
\keyword{links}
\keyword{homology 3--sphere}
\keyword{maximal Euler characteristic}
\keyword{taut foliations}
\keyword{alternative links}
\subject{primary}{msc2000}{57R58}
\subject{primary}{msc2000}{53D40}
\subject{secondary}{msc2000}{57M27}
\subject{secondary}{msc2000}{57R30}

\received{11 June 2005}
\revised{6 January 2006}
\published{21 June 2006}
\publishedonline{21 June 2006}
\accepted{9 May 2006}
\proposed{Peter Ozsv\'ath}
\seconded{Rob Kirby, Jim Bryan}
\corresponding{}
\editor{}
\version{}

\AtBeginDocument{\let\bar\wbar}
\def\co{\mskip0.2mu\colon\mskip2.8mu minus0.8mu\relax}
\newenvironment{CD}{\diagram[amstex]}{\enddiagram}
\diagramstyle[labelstyle=\scriptstyle,heads=curlyvee]
\makeop{Gr}
\makeop{Fil}
\let\mathscr\mathcal

\makeatletter
\def\cnewtheorem#1[#2]#3{\newtheorem{#1}{#3}[section]
\expandafter\let\csname c@#1\endcsname\c@thm}

\newcommand{\homo}{{\mathrm H}}
\newcommand{\dif}{{\mathrm d}}

\newtheorem{thm}{Theorem}[section]
\cnewtheorem{lem}[thm]{Lemma}
\cnewtheorem{cor}[thm]{Corollary}
\cnewtheorem{prop}[thm]{Proposition}

\theoremstyle{definition}
\cnewtheorem{defn}[thm]{Definition}
\cnewtheorem{rem}[thm]{Remark}
\cnewtheorem{algorithm}[thm]{Algorithm}
\cnewtheorem{operation}[thm]{Operation}
\cnewtheorem{notation}[thm]{Notation}

\makeatother  %  move after \newtheorem block

\begin{document}

\begin{asciiabstract}
Ozsvath and Szabo proved that knot Floer homology determines
the genera of knots in S^3. We will generalize this deep result
to links in homology 3-spheres, by adapting their method. Our
proof relies on a result of Gabai and some constructions related
to foliations. We also interpret a theorem of Kauffman in the
world of knot Floer homology, hence we can compute the top
filtration term of the knot Floer homology for alternative links.
\end{asciiabstract}

\begin{htmlabstract}
Ozsv&aacute;th and Szab&oacute proved that knot Floer homology determines
the genera of knots in S<sup>3</sup>. We will generalize this deep result
to links in homology 3&ndash;spheres, by adapting their method. Our
proof relies on a result of Gabai and some constructions related
to foliations. We also interpret a theorem of Kauffman in the
world of knot Floer homology, hence we can compute the top
filtration term of the knot Floer homology for alternative links.
\end{htmlabstract}

\begin{abstract}
Ozsv\'ath and Szab\'o proved that knot Floer homology determines
the genera of knots in $S^3$. We will generalize this deep result
to links in homology 3--spheres, by adapting their method. Our
proof relies on a result of Gabai and some constructions related
to foliations. We also interpret a theorem of Kauffman in the
world of knot Floer homology, hence we can compute the top
filtration term of the knot Floer homology for alternative links.
\end{abstract}

\maketitle

\section{Introduction}\label{sec:1}

The central problem in knot theory is the classification of knots
and links. In order to study this problem, people have introduced
a lot of invariants. Among them, knot Floer homology seems to be
particularly interesting. This invariant was introduced by
Ozsv\'ath and Szab\'o in \cite{OSz4}, and independently by
Rasmussen in \cite{Ra}, as part of Ozsv\'ath and Szab\'o's
Heegaard Floer theory.

One remarkable feature of knot Floer homology is that it detects the
genus in the case of classical knots, Ozsv\'ath and Szab\'o
\cite[Theorem 1.2]{OSz7}, namely, the genus of a classical knot is the
highest nontrivial filtration level of the knot Floer homology. A
direct corollary is: knot Floer homology detects the unknot. Even this
corollary has its own interest.

In this paper, we will generalize this deep result to links in
homology 3--spheres. Our main theorem is the following.

\begin{thm}\label{EulerBound} Suppose $L$ is an oriented link in a closed 3--manifold $Z$ with $\homo_1(Z)=0$.
$|L|$ denotes the number of components of $L$, and $\chi(L)$
denotes the maximal Euler characteristic of the Seifert surface
bounded by $L$. Then
$$\frac{|L|-\chi(L)}2=\max\{\,i\,|\,\widehat{HFK}(Z,L,i)\ne0\,\}.$$
In particular, knot Floer homology detects trivial links in
homology 3--spheres.
\end{thm}

\begin{rem}
When $Z=S^3$, the Euler characteristic of $\widehat{HFK}(S^3,L)$
gives rise to the Alexander--Conway polynomial $\Delta_L(t)$ (see
Ozsv\'ath--Szab\'o \cite{OSz4}). But knot Floer homology certainly contains more
information than $\Delta_L(t)$. For example, when $L$ is a
boundary link, $\Delta_L$ is 0, while our theorem implies
$\widehat{HFK}(S^3,L)$ is nontrivial.
\end{rem}

Our paper is organized as follows:

In \fullref{sec:2} we give some background on Heegaard Floer theory,
especially some useful remarks about links.

In \fullref{sec:3} we will state a result of Gabai in \cite{G2}, which
allows us to prove the main theorem for knots in homology spheres
immediately. But in order to study links, we need some
constructions related to foliations. As a consequence of the
constructions, we can embed a certain 3--manifold into a
symplectic 4--manifold with a ``good" symplectic structure.
Results in this section have their own interest.

\fullref{sec:4} will be devoted to the proof of our main theorem for
links. Our argument is not much different from Ozsv\'ath and
Szab\'o's original one.

In \fullref{sec:5}, we interpret a theorem of Kauffman on alternative
links, using the language of knot Floer homology. Hence we can
compute the top filtration term of the knot Floer homology for
alternative links. This section is independent of our main result.

\medskip
{\bf Acknowledgements}\qua We are grateful to David Gabai,
Peter Kronheimer and Zolt\'an Szab\'o for some helpful
conversations.

The author is partially supported by the Centennial fellowship of
the Graduate School at Princeton University.

\section{Preliminaries}\label{sec:2}

\subsection{Heegaard Floer theory}

We will briefly include some background on Heegaard Floer theory
here, a survey of Heegaard Floer theory can be found in
Ozsv\'ath--Szab\'o \cite{OSz9}.

For any closed oriented 3--manifold $Y$, and a Spin$^c$ structure
$\mathfrak s$ on $Y$, Ozsv\'ath and Szab\'o associate a package of
Floer homologies: $HF^{\infty}(Y,\mathfrak s)$, $HF^+(Y,\mathfrak
s)$, $HF^-(Y,\mathfrak s)$, $\widehat{HF}(Y,\mathfrak s)$,
$HF_{\textrm{red}}(Y,\mathfrak s).$ These are topological
invariants of $(Y,\mathfrak s)$. Furthermore, the homologies
(except $\widehat{HF}(Y,\mathfrak s)$) admit an action $U$, which
lowers degree by $2$. So these groups are $\mathbb Z[U]$--modules.
When $\mathfrak s$ is torsion, one can equip an absolute $\mathbb
Q$--grading to these groups.

This theory is functorial. Namely, if we have a cobordism $W$ from
$Y_1$ to $Y_2$, $\mathfrak r$ is a Spin$^c$ structure on $W$,
$\mathfrak s_i$ is the restriction of $\mathfrak r$ on $Y_i$, then
there is a homomorphism
\begin{equation}F^{\circ}_{W,\:\mathfrak r}\co HF^{\circ}(Y_1\:,\mathfrak s_1)\to HF^{\circ}(Y_2\:,\mathfrak s_2).\label{eq:1}\end{equation}
Here $HF^{\circ}$ denotes any one of the five homologies. When
$\mathfrak s_1,\mathfrak s_2$ are both torsion,
$F^{\circ}_{W,\mathfrak r}$ shifts degree by
\begin{equation}\frac{c_1(\mathfrak r)^2-2\chi(W)-3\sigma(W)}4.
\label{eq:2}\end{equation}
Moreover, suppose $W_1$ is a cobordism from $Y_1$ to $Y_2$, $W_2$
is a cobordism from $Y_2$ to $Y_3$, $W=W_1\cup_{Y_2}W_2$,
$\mathfrak r_i$ is a Spin$^c$ structure on $W_i$. Then we have the
composition law:
\begin{equation}
F^{\circ}_{W_2,\:\mathfrak r_2}\circ F^{\circ}_{W_1,\:\mathfrak r_1}=
\sum_{\mathfrak r\,\in\,\textrm{Spin}^c(W),\:\mathfrak r|W_i=\mathfrak
r_i}\pm F^{\circ}_{W,\:\mathfrak r}
\label{eq:3}\end{equation}
The above $F^{\circ}_{W,\:\mathfrak r}$ is an invariant of
cobordisms. In particular, when $X$ is a closed 4--manifold with
$b_2^+>1$, $\mathfrak r\in\textrm{Spin}^c(X)$, Ozsv\'ath and
Szab\'o can define a mixed invariant $\Phi_{X,\:\mathfrak
r}\in\mathbb Z$. $\Phi$ behaves like, and is conjectured to agree
with, the Seiberg--Witten invariant.

When $K$ is a null-homologous knot in $Y$, one can define a
$\mathbb Z\oplus\mathbb Z$--filtered chain complex
$CFK^{\infty}(Y,K)$. From this chain complex, we can get some knot
invariants, one of which is called the knot Floer homology,
denoted by $\widehat{HFK}(Y,K)$. Given a Seifert surface of $K$,
one can associate a filtration to $\widehat{HFK}(Y,K)$. The
summand of $\widehat{HFK}(Y,K)$ at filtration level $i$ is denoted
by $\widehat{HFK}(Y,K,i)$.

\subsection{A few remarks about links}

In \cite{OSz4}, Ozsv\'ath and Szab\'o gave a well-defined
correspondence from links to knots, and the invariant for links is
defined to be the invariant for the corresponding knots.

The construction is as follows: given a null-homologous oriented
$n$--component link $L$ in $Y$, choose two points $p,q$ on
different components of $L$. Remove two balls at $p,q$, then glue
in a 3--dimensional tube $S^2\times I$. (In \cite{OSz4}, this
$S^2\times I$ is called a 1--handle, but we would rather call it a
tube.) Inside the tube, there is a band, along which we can
perform a connected sum of the two components of L containing $p$
and $q$. We choose the band so that the connected sum respects the
original orientation on $L$. Now we have a link in $Y^3\#S^2\times
S^1$, with one fewer components. Repeat this construction until we
get a knot. The new knot is denoted by $\kappa(L)$, and the new
manifold is denoted by $\kappa(Y)=Y\#(n-1)(S^2\times S^1)$.
Ozsv\'ath and Szab\'o proved that this correspondence
$(Y,L)\mapsto(\kappa(Y),\kappa(L))$ is well-defined.

At first sight, the above construction seems rather artificial.
But as we will see, it is pretty reasonable. As for now, we are
satisfied to justify the construction by giving two folklore
propositions.

\begin{prop}[Adjunction Inequality]\label{adjunct}
$L\subset Y$ is an oriented link. If $\widehat{HFK}(Y,L,i)$ $\ne0$,
then, for each Seifert surface $F$ for $L$, we have that
$$i\le\frac{|L|-\chi(F)}2.$$
\end{prop}
\begin{proof}
Let $\kappa(F)$ be the Seifert surface of the knot $\kappa(L)$,
obtained by adding $|L|-1$ bands to $F$. Then
$\chi(\kappa(F))=\chi(F)-(|L|-1)$, hence
$\frac{|L|-\chi(F)}2=g(\kappa(F))$. Our result holds by the usual
adjunction inequality for knots.
\end{proof}

\begin{prop}\label{fibred}
Let $L\subset Y$ be an oriented fibred link, ie, $Y-L$ fibers
over the circle, and the fiber is a Seifert surface $F$ of $L$.
Then
$$\widehat{HFK}(Y,L,\frac{|L|-\chi(F)}2)\cong\mathbb Z.$$
\end{prop}
\begin{proof}
Suppose $p,q\in L$ are two points in different components. Add a
tube $R=S^2\times I$ to $Y$ with feet at $p,q$, and perform the
band-connected sum operation to $L$, we get a new link $L'\subset
Y\#S^2\times S^1$. Its Seifert surface $F'$ is obtained by adding
a band $B$ to $F$.

$F$ is connected, hence we can find an arc $\gamma\subset F$
connecting $p$ to $q$. A neighborhood of $\gamma$ in $F$ is a
rectangle $D$, we choose $D$ so that $B\cup D$ is an annulus $A$.
$N=\mathrm{Nd}(R\cup\gamma)$ is a punctured $S^2\times S^1$, we
can isotope $F'$ slightly so that $F'\cap\partial N=D$, $F'\cap
N=A$. Cap off the boundary sphere of $N$ by a ball, the new
manifold is denoted by $\widehat N$. Now $\partial A$ is a fibred
link in $\widehat N$, with fiber $A$, and $(Y\#S^2\times S^1, L')$
is obtained by plumbing the two fibred links $(Y,L)$, $(\widehat
N,
\partial A)$. It is well-known that $(Y\#S^2\times S^1, L')$ is
also fibred. This fact was first proved by Stallings in \cite{St},
using algebraic method. A geometric proof can be found in
Gabai \cite{G1}, where it is also shown that if $L'$ is fibred then so
is $L$.

The above proof shows that $(\kappa(Y),\kappa(L))$ is a fibred
knot, now we use Ozsv\'ath--Szab\'o \cite[Theorem 1.1]{OSz5}.
\end{proof}

\subsection{Kauffman states for links}

The material in this subsection will be used in \fullref{sec:5}, not in
the proof of our main theorem.

There is no known algorithm to compute the invariants in Heegaard
Floer theory. But one can find all the generators of the chain
complex in a purely combinatorial way. In \cite{OSz6}, Ozsv\'ath
and Szab\'o described such a method for classical knots.

Suppose $L\subset S^3$ is an oriented link. Choose a generic
projection of $L$, we get a 4--valent graph $\mathcal D$ on $S^2$.
$\mathcal D$ is called a {\it diagram} for $L$. We always choose
$\mathcal D$ so that it is {\sl connected}. If we fix an edge in
$\mathcal D$, then $\mathcal D$ is called a {\it decorated
diagram}.

Suppose a decorated diagram $\mathcal D$ has $m$ vertices, then
$\mathcal D$ divides $S^2$ into $m+2$ regions. Suppose $A,B$ are
the two regions sharing the distinguished edge. A {\it Kauffman
state} is a correspondence $x$, which associates one of the four
corners at $v$ to each vertex $v$ of $\mathcal D$, so that the
corners are not in $A,B$, and any two corners are not in the same
region.

\begin{figure}[ht!]\small
\labellist
\pinlabel 0 at 0 102
\pinlabel 0 at 147 102
\pinlabel 0 at 235 102
\pinlabel 0 at 382 102
\pinlabel $-\frac12$ at 75 203
\pinlabel $\frac12$ at 75 0
\pinlabel $\frac12$ at 309 203
\pinlabel $-\frac12$ at 309 0
\endlabellist
\cl{\includegraphics[width=2in]{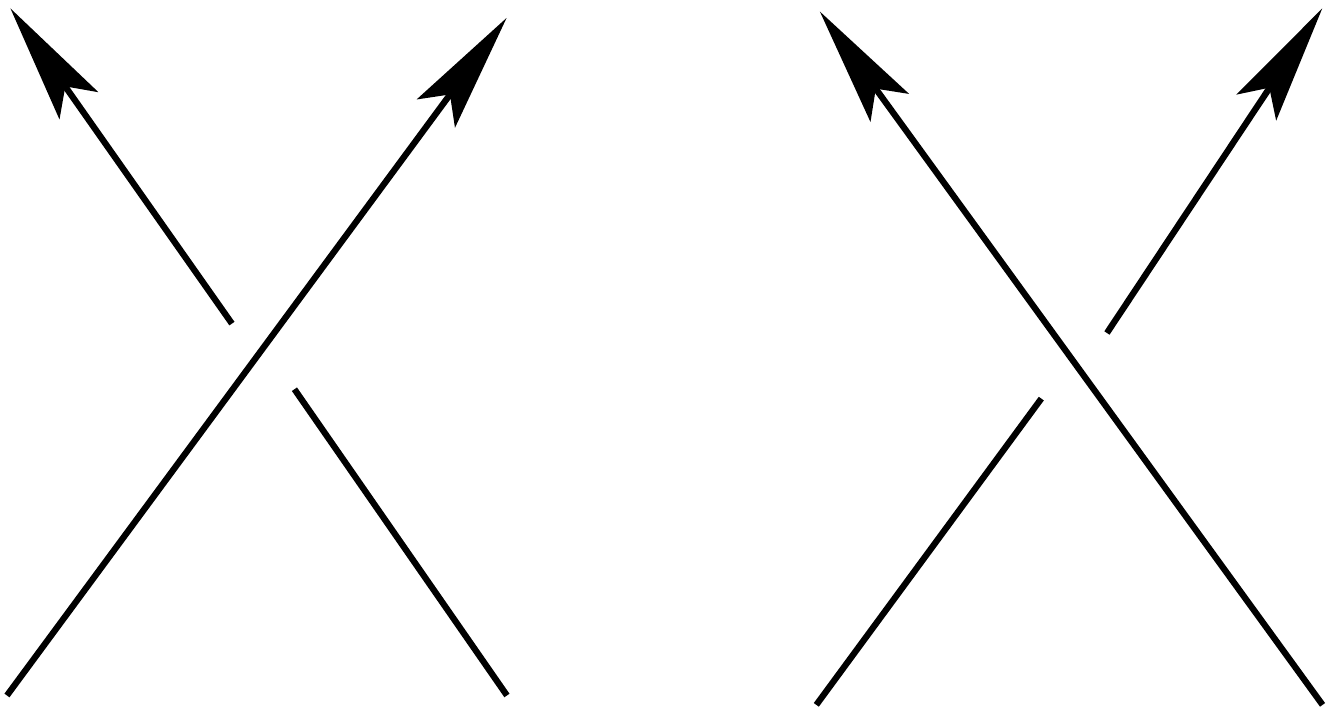}}
\caption{Local filtration level contribution $\Fil(x(v))$}
\label{fig:1}
\end{figure}

One can also associate filtration level and grading to a Kauffman
state. The local filtration level contribution $\Fil(x(v))$ at $v$
is shown in \fullref{fig:1}, and the local grading contribution
$\Gr(x(v))$ at $v$ is shown in \fullref{fig:2}. We define
$$\Fil(x)=\frac{|L|-1}2+\sum_{v\in\textrm{Vert}(\mathcal D)}\Fil(x(v)),$$
$$\Gr(x)=\frac{|L|-1}2+\sum_{v\in\textrm{Vert}(\mathcal D)}\Gr(x(v)).$$
When $L$ is a knot, one can choose the chain complex
$\widehat{CFK}(S^3,L)$, so that the generators are precisely the
Kauffman states with filtration level $\Fil(x)$ and grading $\Gr(x)$
\cite{OSz6}.

\begin{figure}[ht!]\small
\labellist
\pinlabel 0 at 0 102
\pinlabel 0 at 147 102
\pinlabel 0 at 235 102
\pinlabel 0 at 382 102
\pinlabel $-1$ at 75 203
\pinlabel 0 at 75 0
\pinlabel 1 at 309 203
\pinlabel 0 at 309 0
\endlabellist
\cl{\includegraphics[width=2in]{\figdir/cross}}
\caption{Local grading contribution $\Gr(x(v))$}
\label{fig:2}
\end{figure}

The reader should note that, when $L$ is a link, the Kauffman
states do not correspond to the generators of
$\widehat{CFK}(S^3,L)$ in a natural way. So the notions
``filtration level" and ``grading" do not make sense for links. In
\fullref{sec:5}, we will show that, for alternative links, the Kauffman
states at the top filtration level are in one-to-one
correspondence with the generators of the corresponding knot Floer
homology group, and the filtration level and grading of the
Kauffman states are equal to the filtration level and grading of
 the knot Floer homology at the top filtration level.

\begin{rem}
During the course of this work, Khovanov and Rozansky posted a
paper \cite{KhR}, in which they associate a triply-graded
cohomology to a classical link. Euler number of this cohomology
gives the HOMFLYPT polynomial. It is not known whether their work
can eventually lead to a purely combinatorial account of knot
Floer homology for classical links.
\end{rem}

\section{A few constructions related to taut foliations}\label{sec:3}

Suppose $M$ is a 3--manifold, $S\subset M$ is an embedded sphere.
We say $S$ is an {\it essential sphere} if $S$ does not bound a
ball in $M$. $M$ is {\it irreducible} if it does not contain any
essential sphere. Otherwise we say $M$ is {\it reducible}. We will
need the following lemma, whose proof only involves standard
arguments in 3--dimensional topology.

\begin{lem}\label{nonsplit}
$Z$ is a closed oriented 3--manifold, $L\subset Z$ is an
$n$--component null-homologous oriented link,
$(Y,K)=(\kappa(Z),\kappa(L))$. If $Z-L$ is irreducible, then $Y-K$
is irreducible.
\end{lem}
\begin{proof}
Suppose $G$ is a Seifert surface of $K$ with minimal genus.
$R_1,\dots,R_{n-1}$ are the tubes added to $Z$,
$P_1,\dots,P_{n-1}$ are their belt spheres, ie,
$S^2\times\textsl{point}$ in the tube $S^2\times I$. Assume $S$ is
an essential sphere in $Y-K$.

We can assume $P_k,G,S$ are mutually transverse. $P_k\cap K$
consists of two points, hence $P_k\cap G$ consists of an arc and
some circles. If $C$ is an innermost circle in $P_k-\textsl{arc}$,
then $C$ bounds a disk $D$. Cut $G$ open along $C$, glue in two
copies of $D$, and possibly throw away a closed component, we get
a new Seifert surface $G'$. $G'\cap P_k$ has fewer components than
$G\cap P_k$, and $\chi(G')\ge\chi(G)$. Hence we can assume
$P_k\cap G$ is an arc.

$S\cap G$ consists of circles, each circle bounds a disk in $G$
since $G$ is a Seifert surface with minimal genus. Now assume $C$
is an innermost circle in such a disk, $C$ bounds $D$. Using $D$
to do surgery to $S$, we get two spheres $S_1,S_2$. $S$ is
essential in $Y-K$, hence one of $S_1,S_2$ must be essential, and
the intersection of this sphere with $G$ has fewer components.
Hence we can assume $S$ is disjoint from $G$.

By the same argument, we can assume $S$ is disjoint from all
$P_k$'s, hence $S$ lies in $Y-\cup_{k=1}^{n-1}R_k$, which is a
submanifold of $Z$. Hence $S\subset Z-L$. Since $Z-L$ is
irreducible, $S$ bounds a 3--ball, a contradiction.
\end{proof}

\begin{rem}\label{addband}
Suppose $F$ is a Seifert surface of $L$, with maximal Euler
number. The second paragraph in the above proof shows that, a
minimal genus Seifert surface of $K$ can be obtained by adding
$n-1$ bands to $F$.
\end{rem}

In \cite{OSz7}, Ozsv\'ath and Szab\'o proved that, for classical
knots, genus is the highest nontrivial filtration level of the
knot Floer homology. Their proof relies on the existence of a
``nice" taut foliation of knot complement, which was proved by
Gabai in \cite{G2}. In the same paper, Gabai also gave another
existence result for taut foliations:

\begin{thm}\label{foliation}{\rm\cite[Theorem 8.9]{G2}}\qua $K$ is a
null-homologous knot in $Y$. $Y$ is reducible and $Y-K$ is
irreducible. $\homo_1(Y;\mathbb Z)$ is torsion free. If $G$ is a
minimal genus Seifert surface of $K$, then there exists a taut
finite depth foliation $\mathscr F$ of
$Y-\mathrm{int}(\mathrm{Nd}(K))$, such that $\mathscr F|\partial
\mathrm{Nd}(K)$ is a foliation by circles which are longitudes,
and $G$ is a leaf of $\mathscr F$. Moreover, when $G$ has genus
$>1$, one can arrange $\mathscr F$ to be smooth.
\end{thm}

\begin{rem}
We briefly describe how Gabai can extend his result for classical
knots to this case. In $S^3$, there is a Heegaard sphere. Gabai
introduced the notion of ``thin position", hence he could use the
Heegaard sphere like an essential sphere. When $Y$ is reducible,
there is already given an essential sphere, so Gabai's argument
can proceed.

``$\homo_1(Y)$ is torsion free" is a technical condition, to
eliminate some bad cases like the existence of a Scharlemann
cycle.

Gabai did not mention smoothness in the statement of \cite[Theorem
8.9]{G2}; however, one can obtain the smoothness conclusion when
genus $>1$ just as in the proof of the main result of that paper.
\end{rem}

Having \fullref{foliation}, we can immediately prove our main
theorem in the case of knots, namely, the following result.

\begin{prop}\label{knotbound}
Suppose $K$ is a knot in a closed 3--manifold $Z$, $\homo_1(Z)=0$.
Then
$$g(K)=\max\{\,i\in\mathbb Z\,|\,\widehat{HFK}(Z,K,i)\ne0\,\}.$$
\end{prop}
\begin{proof}
If $Z-K$ is reducible, then $Z-K=(Z'-K)\#Z''$, where $Z'-K$ is
irreducible, and $Z',Z''$ are homology spheres. Now we have
$$\widehat{HFK}(Z,K)\cong\widehat{HFK}(Z',K)\otimes\widehat{HF}(Z'').$$
$\widehat{HF}(Z'')$ is nontrivial since its Euler characteristic
is $1$, so we reduce our problem to $(Z',K)$.

From now on, we assume $Z-K$ is irreducible, and $Z$ is not $S^3$.
Consider the knot $K\#K$ in $Z\#Z$. We find that the conditions in
\fullref{foliation} hold, namely, $Y=Z\#Z$ is reducible, the
complement of $K\#K$ is irreducible, and $\homo_1(Y)$ is torsion
free. Hence there exists a taut smooth foliation $\mathscr F$ of
$Y-\mathrm{int}(\mathrm{Nd}(K\#K))$, such that $\mathscr
F|\partial \mathrm{Nd}(K\#K)$ is a foliation by longitudes.

Now we can prove our proposition for $K\#K$, the proof is exactly
same as the proof of \cite[Theorem 1.2]{OSz7}, modulo the
following \fullref{0surgery1}. The reader can find the argument
in \cite{OSz7}, or a more delicate but essentially same argument
in the proof of our main theorem.

By the K\"unneth formula for connected sums \cite[Theorem
7.1]{OSz4}, we get the desired result for $K$.
\end{proof}

\begin{lem}[Compare {\cite[Corollary 4.5]{OSz4}}]\label{0surgery1} 
Let $K$ be a knot in a homology 3--sphere $Y$,
$Y_p$ be the manifold obtained by $p$--surgery on $K$. Let $d$ be
an integer satisfying $\widehat{HFK}(Y,K,i)=0$ when $i\ge d$, and
suppose that $d>1$. Then
$$HF^+(Y_0,d-1)=0.$$
\end{lem}
\begin{proof}
Consider the bigraded chain complex $C=CFK^{\infty}(Y,K)$. We have
the short exact sequence
$$\begin{CD}
0 \to  C\{i<0 \;\textrm{and}\; j\ge d-1\} @>>> C\{i\ge0
\;\textrm{or}\; j\ge d-1\} @>\Psi>>C\{i\ge0\} \to 0,
\end{CD}
$$
hence a long exact sequence relating their homologies. In $C\{i<0
\;\textrm{and}\; j\ge d-1\}$, all summands lie at filtration
levels $\ge d$. By our assumption,
$$\homo_*(C\{i<0 \;\textrm{and}\; j\ge d-1\})=0.$$
Hence $\Psi$ induces an isomorphism $\psi$ on the level of
homology.

By \cite[Theorem 4.4]{OSz4} (also compare \cite[Theorem
2.3]{OSz8}), when $p>0$ is a sufficiently large integer,
$$\homo_*(C\{i\ge0
\;\textrm{or}\; j\ge d-1\})\cong HF^+(Y_p,d-1),$$ and the map
$\Psi$ coincides with the map induced by the cobordism $W$ from
$Y_p$ to $Y$, endowed with the Spin$^c$ structure $\mathfrak r$
with
$$\langle c_1(\mathfrak r), [\widetilde F]\rangle=2d-2-p.$$ Here $F$
is a Seifert surface of $K$, and $\widetilde F$ is a closed
surface in $W$, obtained from $F$ by capping off $\partial F$ in
$W$.

According to \cite[Theorem 9.19]{OSz1}, there is an long exact
sequence
$$\begin{CD}
\cdots\to HF^+(Y_0\:,d-1) @>>>
HF^+(Y_p\:,[d-1])@>f>>HF^+(Y) \to\cdots
\end{CD}$$
We will compare $\psi$ and $f$.

Suppose $\mathfrak r'\in \textrm{Spin}^c(W)$ nontrivially
contributes to $f$, then $\mathfrak r'-\mathfrak r$ must be
$k\textrm{PD}([\widetilde F])$ for some integer $k$, where
$\textrm{PD}$ is the Poincar\'e duality map. It is not hard to
find
$$c_1^2(\mathfrak r')=-\frac{(2d-2-p+2kp)^2}p,\quad k\in\mathbb Z.$$
Now by the degree shifting formula \ref{eq:2}, $f$ has the form
$\psi+\iota$, where $\iota$ is a sum of homogeneous maps which
have lower orders than $\psi$. Since $\psi$ is an isomorphism, it
is clear that $\psi+\iota$ is also an isomorphism, hence
$HF^+(Y_0\:,d-1)=0$.
\end{proof}

\begin{rem}
Now we can prove our main theorem for links which have connected
minimal Seifert surfaces, by a pretty trivial argument sketched as
follows. If $F$ is a connected minimal Seifert surface of
$L\subset Z$, then we can perform Murasugi sums of $F$ with Hopf
bands, to get a knot $K$ in $Z$. Using the skein exact sequence
and the adjunction inequality, we can reduce our problem to $K$,
for which the problem is already solved in \fullref{knotbound}.
\end{rem}

In order to prove our main theorem for general links, we will do
some constructions related to foliations, contact structures and
symplectic structures.

\begin{rem}\label{ETElEt} We recall the procedures in \cite{El,ElT,Et} 
of passing from taut foliations to symplectic
structures. Suppose $\mathscr G$ is a $C^2$ taut foliation of $M$,
then there exists a closed 2--form $\omega$ on $M$, so that
$\omega$ does not vanish on $\mathscr G$. Let $W=M\times I$,
define a closed 2--form
${\Omega'}=p^*\omega+\varepsilon\dif(t\alpha)$. Here $p$ is the
projection of $W$ onto $M$, $\alpha$ is the 1--form defining
$\mathscr G$. Hence $(W,\Omega')$ weakly symplectically semi-fills
$(M,\mathscr G)$. The plane field $\mathscr G$ can be slightly
perturbed to a contact structure, which is also weakly
symplectically semi-fillable.

By \cite[Theorem 1.3]{El}, $(W,\Omega')$ can be symplectically
embedded into a closed symplectic manifold $(X,\Omega)$.
\end{rem}

We also need the following basic operation of modifying
foliations, which preserves tautness.

\begin{operation}[Creating or killing holonomy]
\label{creatholo}
Suppose $\mathscr G$ is a foliation of $M$, $S^1\times I$ is an
annulus in $M$, so that $\mathscr G$ restricts to $S^1\times t$'s
on $S^1\times I$. $f\co I\to I$ is a homeomorphism supported in
$\mathrm{int}(I)$. We can cut open $M$ along $S^1\times I$, then
reglue by $\mathrm{id}\times f$, thus get a new foliation
$\mathscr G'$. If $\mathscr G$ and $f$ are smooth, then $\mathscr
G'$ is also smooth.
\end{operation}

In order to prove our main theorem, it seems that \fullref{foliation} and \fullref{ETElEt} are not enough. We will
also need a ``good" symplectic structure. The next lemma gives a
way to construct such symplectic structures.

\begin{lem}\label{symplemma}
Let $\mathscr G$ be a smooth taut foliation of $M$.
$T_1,\dots,T_{m}$ are closed surfaces in $M$. Suppose there exist
simple closed curves $c_1,\dots,c_m$, each curve is contained in a
leaf of $\mathscr G$, and
$$[c_i]\cdot [T_j]=\delta_{ij}.$$
If for each curve $c_i$, the holonomy along $c_i$ is trivial, then
$M$ can be embedded into a closed symplectic 4--manifold
$(X,\Omega)$, so that $X=X_1\cup_M X_2$, $b_2^+(X_j)>0$, and
$$\int_{T_i}\Omega=0.$$
\end{lem}
\begin{proof}
Since the holonomy along $c_i$ is trivial, we can choose a
neighborhood $\mathrm{Nd}(c_i)=S^1\times [-1,1]\times[-1,1]$, with
coordinates $(x,y,t)$, so that restriction of $\mathscr G$ to
$\mathrm{Nd}(c_i)$ consists of $S^1\times [-1,1]\times t$'s.

Let $u$ be a bump function supported in $(-1,1)$, so that
$\int_{-1}^1 u(t)\dif t=1$. Hence
$$\varphi_i=u(y)u(t)\dif y\wedge\dif t$$
is a closed 2--form supported in $\mathrm{Nd}(c_i)$. We can extend
$\varphi_i$ by 0 to the whole $M$, then $\varphi_i$ restricts to 0
on $\mathscr G$, and $\int_T \varphi_i=[c_i]\cdot [T]$ for any
closed surface $T\subset M$.

Now we perturb $\omega$ by multiples of
$\varphi_1,\dots,\varphi_m$, so as to get a closed 2--form
$\omega'$, which never vanishes on $\mathscr G$, and satisfies
$\int_{T_i}\omega'=0$ for any $i$.

As in \fullref{ETElEt}, the form
$${\Omega'}=p^*\omega'+\varepsilon\dif(t\alpha)$$
is a symplectic form on $M\times I$. Moreover,
$\int_{T_i}\Omega'=0$ for any $i$. Hence the result holds by
Eliashberg and Etnyre's theorem.
\end{proof}

\begin{notation}
In the rest of this section and \fullref{sec:4}, we will use the
following notation. Let $L$ be an $n$--component link in a
homology sphere $Z$, $n>1$, $(Y,K)=(\kappa(Z),\kappa(L))$. $G$ is
a minimal genus Seifert surface of $K$. By \fullref{addband},
we can assume $G$ is obtained by adding $n-1$ bands
$B_1,\dots,B_{n-1}$ to a Seifert surface $F$ of $L$ with maximal
Euler characteristic. Hence $\chi(G)=\chi(F)-(n-1)$. Let $\widehat
G$ be the extension of $G$ in $Y_0$ obtained by gluing a disk to
$G$.

$Y$ is obtained from $Z$ by adding $n-1$ (3--dimensional) tubes
$R_1,\dots,R_{n-1}$. Suppose $P_i$ is the belt sphere of the tube
$R_i$. The knot $K$ intersects $P_i$ in exactly 2 points, we can
remove two disks from $P_i$ at these two points, then glue in a
long and thin (2--dimensional) tube along an arc in $K$, so as to
get a torus $T_i$. $T_i$ is homologous to $P_i$, but disjoint from
$K$. $\;T_i$ will play an important role in our proof.
\end{notation}

\begin{prop}\label{sumwithfibred}
Notation is as above. We suppose that $Z-L$ is irreducible, and
$\textsl{genus}(G)>1$. After doing connected sum with some fibred
knots in $S^3$, we get a new link $L^*$. We consider
$(Y^*,K^*)=(\kappa(Z),\kappa(L^*))$, and the 0--surgered space
$Y_0^*$. The conclusion is: for a suitably chosen $L^*$, $Y^*_0$
can be embedded into a closed symplectic 4--manifold $(X,\Omega)$,
so that $X=X_1\cup_{Y^*_0}X_2$, $b_2^+(X_j)>0$, and
$$\int_{T_i^*}\Omega=0$$
for all $i$. Moreover, $$ \langle c_1(\mathfrak
k(\Omega)),[\widehat{G^*}]\rangle=2-2g(\widehat{G^*}).$$
\end{prop}
\begin{proof}
We first consider the case that $L$ consists of two components
$K_1,K_2$. We apply \fullref{foliation} to get a smooth taut
foliation $\mathscr F$ on $Y-\mathrm{int}(\mathrm{Nd}(K))$.
According to Gabai \cite{G2,G3}, each leaf of $\mathscr
F$ intersects $\partial \mathrm{Nd}(K)$.

$P_1$ is the belt sphere of the (3--dimensional) tube $R_1$.
$P_1-\mathrm{int}(\mathrm{Nd}(K))$ is an annulus. A standard
argument (see \cite[Theorem 4]{T}) shows that one can deform $P_1$
so that $\mathscr F|(P_1-\mathrm{int}(\mathrm{Nd}(K)))$ is a
product foliation consisting of arcs. Note that
$Y-\mathrm{int}(\mathrm{Nd}(K))-\mathrm{int}(\mathrm{Nd}(P_1))$ is
diffeomorphic to $Z-\mathrm{int}(\mathrm{Nd}(L))$, hence $\mathscr
F$ induces a foliation $\mathscr E=\kappa^{-1}(\mathscr F)$ of
$Z-\mathrm{int}(\mathrm{Nd}(L))$. Conversely, $\mathscr F$ is
obtained by adding bands to leaves of $\mathscr E$.

Each leaf of $\mathscr E$ intersects either $\partial
\mathrm{Nd}(K_1)$ or $\partial \mathrm{Nd}(K_2)$. Let
$$U_i=\{\,x\in Z-\mathrm{int}(\mathrm{Nd}(L))\,|\;\textrm{the leaf containing $x$ intersects $\partial \mathrm{Nd}(K_i)$}\,\}.$$
$U_1,U_2$ are open subsets of $Z-\mathrm{int}(\mathrm{Nd}(L))$,
$U_1\cup U_2=Z-\mathrm{int}(\mathrm{Nd}(L))$. Hence $U_1\cap
U_2\ne\emptyset$, which means that there exists a leaf $\mathcal
L_0$ intersecting both $\partial \mathrm{Nd}(K_1)$ and $\partial
\mathrm{Nd}(K_2)$.

If there is a band $B\subset R_1$, which connects $\mathcal L_0$
to itself, then there is a simple closed curve
$c_1\subset\kappa(\mathcal L_0)$, so that $[c_1]\cdot[T_1]=1$.
Here $\kappa(\mathcal L_0)$ is the leaf of $\mathscr F$ that
contains $\mathcal L_0$.

Now assume that there is a band $B\subset R_1$ connecting
$\mathcal L_0$ to a different leaf $\mathcal L_1$. We can suppose
$B\cap\mathcal L_1\subset\partial \mathrm{Nd}(K_1)$. Perform
connected sum of $L$ and a trefoil $Tr_1\subset S^3$, so that
$Tr_1$ is added to $K_1$. We get a new link $L'$. $\mathscr E$ can
be extended to a foliation $\mathscr E'$ of
$Z-\mathrm{int}(\mathrm{Nd}(L'))$, so that the part of $\mathscr
E'$, obtained from $S^3-Tr_1$, is a fibration over the circle,
with fiber a punctured torus. $\mathcal L_0$ and $\mathcal L_1$
are extended to two leaves $\mathcal L_0'$ and $\mathcal L_1'$ of
$\mathscr E'$. Choose a non-separating circle $c$ on the punctured
torus, extend it to a vertical annulus $c\times I$ in the
fibration part, so that $c\times I$ intersects both $\mathcal
L_0'$ and $\mathcal L_1'$, but $c\times I$ does not intersect the
leaf containing $F$. Then we perform \fullref{creatholo}
along $c\times I$ to $\mathscr E'$, so as to get a new foliation
$\mathscr E''$, where there is a leaf $\mathcal L_0''$ satisfying
 that the band $B\subset R_1$ connects $\mathcal L_0''$ to
itself. See \fullref{fig:3} for an illustration of the local picture.

\begin{figure}[ht!]\small
\labellist 
\pinlabel {$\mathcal L_0'$} [r] at 54 541
\pinlabel {$\mathcal L_1'$} [r] at 54 512
\pinlabel {$\mathcal L_0''$} [r] at 335 541
\endlabellist
\cl{\includegraphics[width=4in]{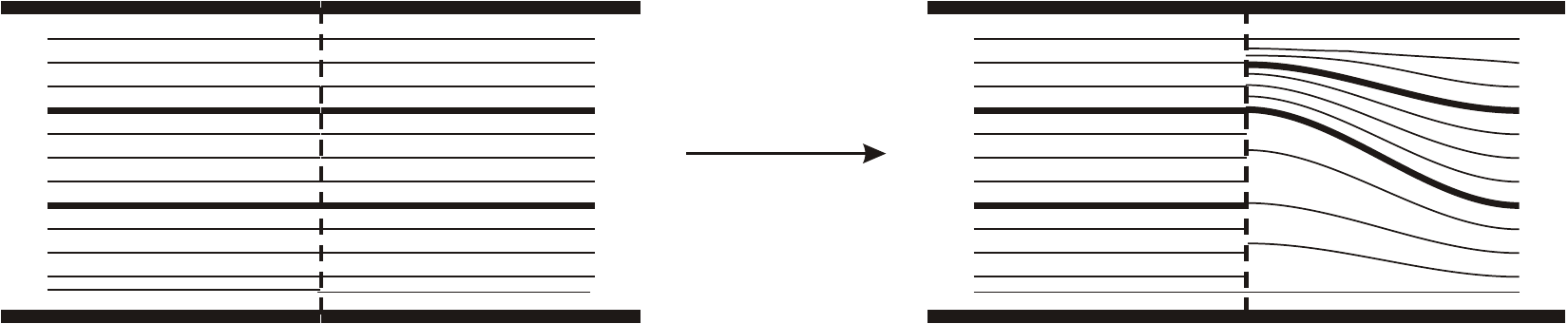}}
\caption{Perform \fullref{creatholo} to merge $\mathcal
L_0'$ and $\mathcal L_1'$}
\label{fig:3}\end{figure}

In any case, we get a simple closed curve
$c_1^{\circ}\subset\kappa(\mathcal L_0^{\circ})$, $c_1^{\circ}$
intersects $T_1^{\circ}$ once. The holonomy along $c_1^{\circ}$
may be nontrivial. We can perform the connected sum of $L^{\circ}$
and a trefoil $Tr_2 \subset S^3$ to get a link $L^*$, and extend
the foliation $\mathscr E^{\circ}$ to a foliation $\mathscr E^*$
of $Z-\mathrm{int}(\mathrm{Nd}(L^*))$ by a fibration over the
circle. Choose two curves $d,e$ on a fiber of the fibration, so
that they intersect in one point. Perform connected sum of
$c_1^{\circ}$ and $d$ to get a new curve $c_1^*$ in a leaf of
$\kappa(\mathscr E^*)$. Do \fullref{creatholo} along
$e\times I$ to cancel the holonomy along $c_1^{\circ}$, then the
holonomy along $c_1^*$ is trivial.

We note that $F^*$ is a leaf or disjoint union of leaves of the
new foliation of $Z-\mathrm{int}(\mathrm{Nd}(L))$, hence $G^*$ is
a leaf of the newly induced foliation of
$Y-\mathrm{int}(\mathrm{Nd}(K^*))$.

Now apply \fullref{symplemma} to get our desired result for the
case that $L$ has two components. The general case can be proved
by the same method. Indeed, by performing connected sum with
fibred knots and creating holonomy, one can get a leaf which
intersects every $\partial \mathrm{Nd}(K_i)$. Again, by the same
procedure one can get a leaf of $\mathscr E$ which is connected to
itself by $n-1$ bands in the $n-1$ different tubes
$R_1,\dots,R_{n-1}$. Hence there are $n-1$ closed curves in the
corresponding leaf of $\kappa(\mathscr E)$, so that they are dual
to $T_1,\dots,T_{n-1}$. Then perform connected sums with fibred
knots and kill the holonomy along these curves. The last step is
to apply \fullref{symplemma}.
\end{proof}

\section{Proof of the main theorem for links}\label{sec:4}

For simplicity, we only prove the case when $L$ is a two-component
link. The reader will find that the proof also works for links
with more components.

Suppose $Y=\kappa(Z)=Z\#S^2\times S^1$, $K=\kappa(L)$. Let
$Y_p=Y_p(K)$ be the manifold obtained from $Y$ by $p$--surgery on
$K$.

The next lemma is an analogue of \fullref{0surgery1}.

\begin{lem}\label{0surgery} Let $d$ be
an integer satisfying $\widehat{HFK}(Y,K,i)=0$ for $i\ge d$, and
suppose that $d>1$. Then
$$HF^+(Y_0,[d-1])=0,$$
where
$$HF^+(Y_0,[d-1])=\bigoplus_{\langle c_1(\mathfrak s),[\widehat G]\rangle=2(d-1)}HF^+(Y_0,\mathfrak s)$$
\end{lem}
\begin{proof}
As in the proof of \fullref{0surgery1}, the chain map
$$\Psi\co C\{i\ge0
\;\textrm{or}\; j\ge d-1\} \to C\{i\ge0\}
$$
induces an isomorphism $\psi$ on the homologies. The map $\Psi$
coincides with the map induced by the cobordism $W$ from $Y_p$ to
$Y$, endowed with the Spin$^c$ structure $\mathfrak r$ with
$$\langle c_1(\mathfrak r), [\widetilde G]\rangle=2d-2-p.$$ Here
$\widetilde G$ is a closed surface in $W$, obtained from $G$ by
capping off $\partial G$ in $W$.

We also have the long exact sequence
$$\begin{CD}
\cdots\to HF^+(Y_0\:,[d-1]) @>>>
HF^+(Y_p\:,[d-1])@>f>>HF^+(Y)\to\cdots
\end{CD}$$
Apply the adjunction inequality \cite[Theorem 7.1]{OSz1} to
$T_1\subset Y_p$, we find that if $\mathfrak t\in
\textrm{Spin}^c(Y_p)$ satisfies $HF^+(Y_p\:,\mathfrak t)\ne0$,
then $\mathfrak t$ is torsion. Hence if $\mathfrak r'\in
\textrm{Spin}^c(W)$ nontrivially contributes to $f$, $\mathfrak
r'-\mathfrak r$ must be a multiple of $\textrm{PD}([\widetilde
G])$, where $\textrm{PD}$ is the Poincar\'e duality map. It is not
hard to find
$$c_1^2(\mathfrak r')=-\frac{(2d-2-p+2kp)^2}p,\quad k\in\mathbb Z.$$
Now by the degree shifting formula \ref{eq:2}, $f$ has the form
$\psi+\iota$, where $\iota$ is a sum of homogeneous maps which
have lower orders than $\psi$. Since $\psi$ is an isomorphism, so
is $\psi+\iota$. Hence $HF^+(Y_0,[d-1])=0$.
\end{proof}

Now we can proceed to the proof of our main theorem. The proof is
not much different from the proof of \cite[Theorem 1.2]{OSz7}.

\begin{proof}[Proof of \fullref{EulerBound}] Suppose $L_1,L_2$
are links in $Z_1,Z_2$, respectively. We have
\begin{eqnarray*}
\widehat{HFK}(Z_1,L_1)\otimes\widehat{HF}(Z_2)&\cong&\widehat{HFK}(Z_1\#Z_2,L_1),\\
\widehat{HFK}(Z_1\#Z_2,L_1\#L_2)\otimes\widehat{HF}(S^2\times
S^1)&\cong&\widehat{HFK}(Z_1\#Z_2,L_1\sqcup L_2).
\end{eqnarray*}
By the above formulas, we can assume $Z-L$ is irreducible. By
doing connected sum with the trefoil, we can assume the genus of
$K$ is bigger than 1. Now apply \fullref{sumwithfibred} to
get a symplectic 4--manifold $(X,\Omega)$, $X=X_1\cup_{Y_0^*}
X_2$, with $b_2^+(X_j)>0$, $\int_{T_1^*}\Omega=0$, and
$$ \langle
c_1(\mathfrak
k(\Omega)),[\widehat{G^*}]\rangle=2-2g(\widehat{G^*}).$$
By the composition formula \ref{eq:3}, the sum
\begin{equation}
\sum_{\eta\in\homo^1(Y_0^*)}\Phi_{X,\mathfrak k(\Omega)+\delta\eta}
\label{eq:4}\end{equation}
is calculated by a homomorphism which factors through
$HF^+(Y_0^*\:,\mathfrak k(\Omega)|_{Y_0^*})$.

$\homo^1(Y_0^*)\cong \mathbb Z\oplus\mathbb Z$ is generated by the
Poincar\'e duals of $[T_1^*]$ and $[\widehat{G^*}]$. So the
Spin$^c$ structures in \ref{eq:4} are precisely
$$\mathfrak
k(\Omega)+a\:\textrm{PD}([T_1^*])+b\:\textrm{PD}([\widehat{G^*}])
\quad(a,b\in\mathbb Z).$$ Here $\textrm{PD}$ is the Poincar\'e
duality map in $X$. The first Chern classes of these Spin$^c$
structures are
$$c_1(\mathfrak
k(\Omega))+2a\:\textrm{PD}([T_1^*])+2b\:\textrm{PD}([\widehat{G^*}]).$$
The evaluation of $c_1(\mathfrak k(\Omega))$ on $[\widehat{G^*}]$
is $2-2g\ne0$. By the degree shifting formula \ref{eq:2}, we conclude
that the terms in \ref{eq:4} that have the same degree as
$\Phi_{X,\mathfrak k(\Omega)}$ are precisely those corresponding
to $\mathfrak k(\Omega)+a\:\textrm{PD}([T_1^*])$. By \cite[Theorem
1.1]{OSz3} and the fact that $\int_{T_1^*}\Omega=0$,
$\Phi_{X,\mathfrak k(\Omega)}$ is the only nontrivial term at this
degree. So $HF^+(Y_0^*\:,\mathfrak k(\Omega)|_{Y_0^*})$ is
nontrivial. Now apply \fullref{0surgery}, we get our desired
result for $L^*$.

The result for $L$ holds by the connected sum formula.\end{proof}

\section{Alternative links}\label{sec:5}

In \cite{Ka}, Kauffman defined a class of links called
``alternative links".  The definition is as follows:

\begin{defn}
Suppose we have a link diagram $\mathcal D$ and that $\mathcal D$ is
connected. We apply Seifert's algorithm (see \cite{Ro}) to
$\mathcal D$, and thus get a collection $\mathcal S$ of disjoint
circles on $S^2$. The vertices of $\mathcal D$ naturally lie in
$S^2-\mathcal S$. We say $\mathcal D$ is an {\it alternative
diagram} if for each component of $S^2-\mathcal S$, the vertices
in this component have the same sign. A link $L$ is called an {\it
alternative link}\; if it has an alternative diagram.
\end{defn}

This definition simultaneously generalizes alternating links and
positive links. Kauffman proved that for an alternative link, the
surface $F$ obtained from Seifert's algorithm has the maximal
Euler characteristic among all the Seifert surfaces. Moreover, the
rank of $\homo_1(F)$ is equal to the degree of the
Alexander--Conway polynomial. Kauffman proved his result by
showing that all the Kauffman states at the top filtration level
have the same grading. The readers are referred to \cite[Chapter
9]{Ka} for a proof. Here, we only state Kauffman's algorithm of
getting all the Kauffman states at the top filtration level.

\begin{algorithm}(Alternative Tree Algorithm, or ATA)

{\bf Step 1}\qua Suppose $\mathcal D$ is a decorated alternative
diagram, choose a checkerboard coloring of the regions. There are
two regions $R_B,R_W$ sharing the distinguished edge, with
coloring black and white respectively. Now we have the black graph
$\Gamma_B$, whose vertices are the regions with black color, and
the edges correspond to vertices of $\mathcal D$. Similarly, there
is the white graph $\Gamma_W$.

{\bf Step 2}\qua For $\Gamma_B$, remove all the edges connecting
two corners with $0$ contribution to the filtration level, see
\fullref{fig:1}. Orient the rest edges so that they point to the corners
with contribution $\frac12$, hence we get an oriented graph
$\Gamma'_B$. $\Gamma'_B$ usually has many components. Similarly,
we have $\Gamma'_W$.

{\bf Step 3}\qua Find a maximal oriented tree $T^1_B\subset
\Gamma'_B$, with root at $R_B$. More precisely, $T^1_B$ is a
maximal tree so that all edges at $R_B$ point away from $R_B$, and
no two edges point to the same vertex. Similarly find a maximal
oriented tree $T^1_W$ with root at $R_W$.

{\bf Step 4}\qua The edges we have chosen correspond to some
vertices of $\mathcal D$. Choose an edge corresponding to one of
the rest vertices, so that it connects a vertex of $T^1_B$ to a
vertex $R^2_B$ in $\Gamma'_B-T^1_B$. Orient this edge so that it
points to $R^2_B$. Now we choose a maximal oriented tree in
$\Gamma'_B$ with root at $R^2_B$. The union of the two black trees
and the edge connecting them is a larger tree, called $T^2_B$.
Similarly, we get a white tree $T^2_W$.

{\bf Step 5}\qua Repeat the process in the last step to enlarge
$T^2_B$ to a tree $T^3_B$, $T^2_W$ to $T^3_W$. Go on with this
process until we get a maximal oriented black tree
$T_B\subset\Gamma_B$ and a maximal oriented white tree
$T_W\subset\Gamma_W$. Suppose $v$ is a vertex of $\mathcal D$, $v$
corresponds to an oriented edge $e$ in $T_B$ or $T_W$. Associate
to $v$ the corner $e$ points to, hence we get a Kauffman state
$x$.
\end{algorithm}

Kauffman proved that there always exist Kauffman states obtained
by ATA, and the Kauffman states so obtained are precisely those
states at the top filtration level. Moreover, they have the same
grading. In fact, we can compute the grading and filtration level
as follows: Suppose $c_-$ is the number of components of
$S^2-\mathcal S$ in which all the vertices of $\mathcal D$ have
negative sign, $r_-$ denotes the number of regions in these
components. Similarly define $c_+,r_+$. Then the grading of the
Kauffman states obtained by ATA is
\begin{equation}
\Gr_{\max}=\frac{|L|-1}2+r_--c_-
\label{eq:5}\end{equation}
and their filtration level is
\begin{equation}
\Fil_{\max}=\frac{|L|-1+r-c}2=\frac{|L|+1+m-c}2
\label{eq:6}\end{equation}
where $r=r_++r_-$, $m=r-2$ is the number of vertices in $\mathcal
D$, $c=c_++c_-$.

To conclude, we interpret Kauffman's work in the world of knot
Floer homology.

\begin{figure}[ht!]\small
\labellist 
\pinlabel {$L_+$} at 207 296
\pinlabel* {$L_-$} [t] at 80 9
\pinlabel* {$L{=}L_0$} [t] at 345 9
\endlabellist
\cl{\includegraphics[width=2in]{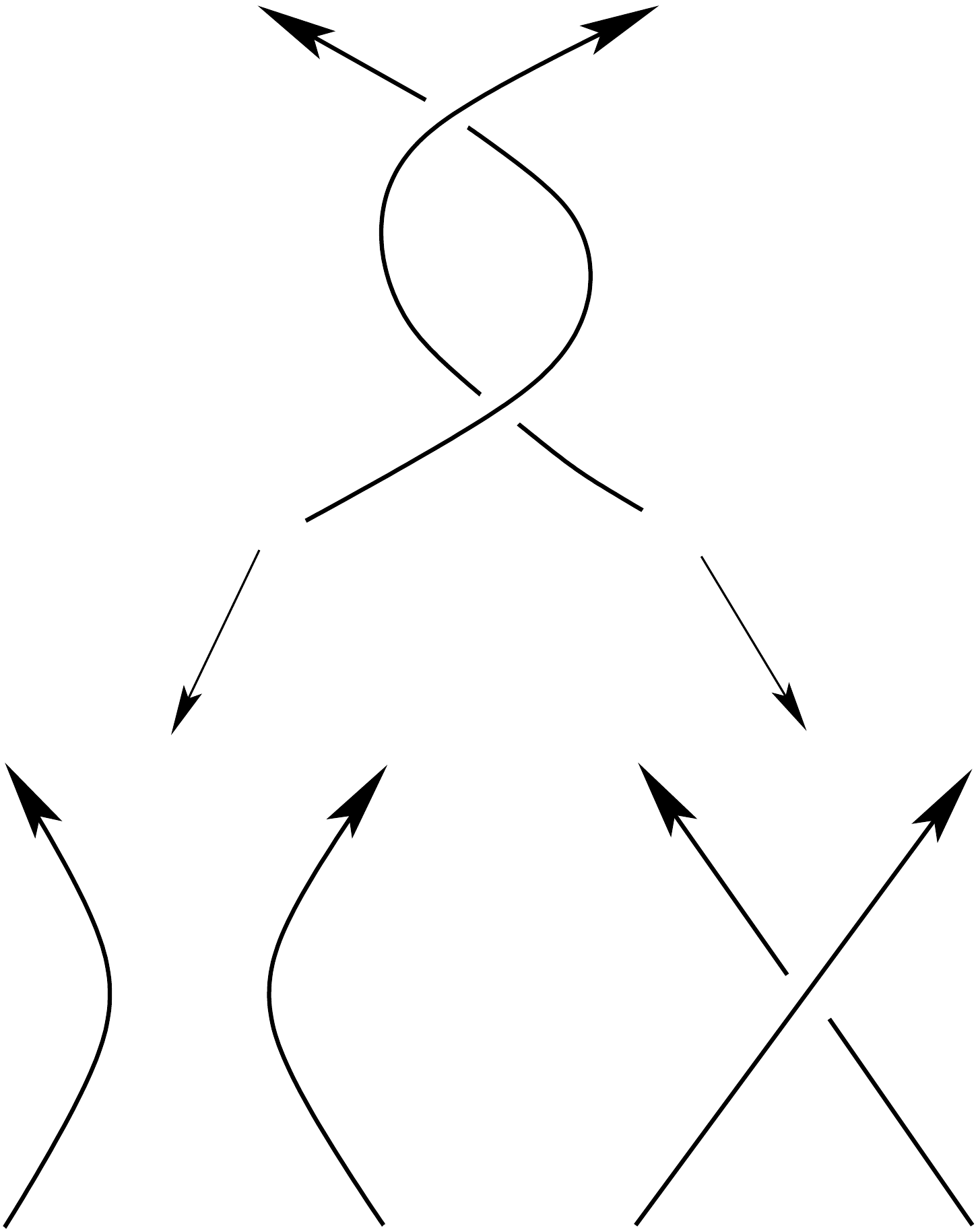}}
\caption{Crossing change}
\label{fig:4}\end{figure}

\begin{thm}
Use notation as above. Suppose $\mathcal D$ is the diagram of an
alternative link $L$, then the maximal nontrivial filtration level
of $\widehat{HFK}(L)$ is the $\Fil_{\max}$ as in (6). Moreover,
$\widehat{HFK}(L,\Fil_{\max})$ is freely generated by the Kauffman
states obtained from ATA, with grading $\Gr_{\max}$ as in (5).
\end{thm}
\begin{proof}
We prove this theorem by induction on $|L|$. When $L$ is a knot,
the result holds by Kauffman's work and the results in Subsection
2.3.

When $L$ is a link, choose two different components which have an
intersection in $\mathcal D$. Without loss of generality, we can
assume this intersection is positive. As in \fullref{fig:4}, $L$ is the
$L_0$, and there are two alternative links $L_+,L_-$ with fewer
components. We have the skein exact sequence
$$
\begin{CD}
\cdots\to\widehat{HFK}(L_-,i)@>>>\widehat{HFK}(L_0,i)@>f>>\widehat{HFK}(L_+,i)
\to\cdots,
\end{CD}
$$
where $f$ drops the grading by $\frac12$.

By induction, $\Fil_{\max}$ is the maximal nontrivial filtration
level of $\widehat{HFK}(L_+)$, and $\widehat{HFK}(L_+,\Fil_{\max})$
is supported at grading $\Gr_{\max}-\frac12$. And $\Fil_{\max}-1$ is
the maximal filtration level of $\widehat{HFK}(L_0)$. Furthermore,
it is easy to find a one-to-one correspondence between the
Kauffman states of $L_+$ and Kauffman states of $L_0$ at the top
filtration level. Now our desired result holds.
\end{proof}

\begin{rem}
In \cite{Cr}, Cromwell proved that an alternative link is fibred
if and only if its Alexander--Conway polynomial is monic, hence if
and only if $\widehat{HFK}(L,\Fil_{\max})\cong\mathbb Z$. Cromwell
actually proved his result for a (possibly but unlikely larger)
class of links called ``homogeneous links".

It is interesting to ask if the converse of \fullref{fibred} is true for 
classical links.
\end{rem}

\bibliographystyle{gtart}
\bibliography{link}

\end{document}